\documentclass[12pt]{amsart}
\usepackage{amsmath}

\newtheorem{Theorem}{Theorem}
\newtheorem{Definition}{Definition}[section]
\newtheorem{Proposition}[Definition]{Proposition}
\newtheorem{Lemma}[Definition]{Lemma}
\newtheorem{Corollary}[Definition]{Corollary}

\newtheorem{Remark}[Definition]{Remark}

\newtheorem{Example}[Definition]{Example}

\newtheorem{Conjecture}[Definition]{Conjecture}

\def\bm{ \left[ \begin{array}{ccc} }
\def\ema{\end{array} \right] }
\def\bmt{ \left[ \begin{array}{cc} }

\newcommand{\R}{{\mathbb R}}
\newcommand{\Z}{{\mathbb Z}}

\newcommand{\C}{{\mathbb C}}

\begin{document}

\begin{center}
\large
Rigidity of locally free Lie group actions\\
 and leafwise cohomology
\end{center}

\vspace{5mm}

\begin{center}
Shigenori Matsumoto (Nihon University, Tokyo)
\end{center}

\vspace{10mm}

\section{Introduction}
In this talk, manifolds, maps, group actions {\em e.\ t.\ c.} are all assumed 
to be of class $C^\infty$. Throughout the talk, $M$ stands for a closed
$C^\infty$ manifold, and $G$ for a connected and simply connected Lie group.
Denote by ${\mathcal A}^r$ the set of locally free right $G$-actions 
(of class $C^\infty$) endowed with the Whitney $C^r$-topology ($1\leq
r\leq \infty$).

An action $\varphi:M\times G \rightarrow M$ in ${\mathcal A}^r$ is
said to be $C^r$ {\em locally rigid} if there exists a neighbourhood
$\mathcal U$ of $\varphi$ in ${\mathcal A}^r$ such that  any $\psi
\in \mathcal U$ is smoothly conjugate to $\varphi$
by a diffeomorphism $F$ of $M$, up to an automorphism $A$ of $G$, {\em i.\ e.}
\begin{equation}  \label{e1}
F(\varphi(x;g))=\psi(F(x);A(g)) 
\end{equation}
for any $x\in M$ and $g\in G$.

An action $\varphi$ is said to be {\em globally rigid} if (\ref{e1})
holds for any $\psi\in{\mathcal A}^r$.

\bigskip

Of course this is an extremely
 strong property. For example any flow (an $\R$-action)
on a manifold $M$ other than $S^1$ cannot be $C^1$ locally rigid.
This follows from Pugh's closing lemma in case the flow does not admit
a periodic orbit. In the other case notice that
the eigenvalues of the Poincar\'e map along
a periodic orbit can be easily changed by a perturbation.

However for other Lie groups, there exist examples of even
globally rigid actions.
Let $GA$ be the Lie group of the orientation preserving affine transformations
on the real line. Let $A\in SL(2;\Z)$ be a hyperbolic automorphism of the
2-torus and denote by ${\mathbb T}_A$  the mapping torus of $A$. 
Then the weak stable foliation of the suspension flow is the orbit foliation
of a locally free $GA$ acion. This action is known to be globally rigid
(\cite{GS,G}).

There are other examples. Let $\Gamma$ be a Fuchsian triangle group.
Then on the manifold $\Gamma\setminus PSL(2;\R)$, a locally free $GA$-action
is defined by the right action of the elements
$$
\bmt
e^{t/2} & x \\ 0 & e^{-t/2}
\ema.
$$ 
The orbit foliation of this action is again a weak stable foliation of
an Anosov flow. E. Ghys (\cite{G})
showed that this action is also globally rigid.

The proofs of these facts consist of the studies of two independent
phenomena. One is concerned about the $C^\infty$ rigidity of the orbit foliation,
and the other is about the rigidity of the parametirization of the action.
When A. Katok and Lewis (\cite{KL}) showed the $C^1$ local rigidity for
 certain $\R^n$-actions, the same strategy was taken. So let us devide the
local rigidity into two parts.

\begin{Definition} {\em An action $\varphi\in{\mathcal A}^r$ is said to be $C^r$
{\em locally orbit rigid} if there exists a neighbourhood $\mathcal U$ of
$\varphi$ such that the orbit foliation ${\mathcal O}_\psi$ of any element
 $\psi\in
\mathcal U$ is smoothly conjugate to the orbit foliation ${\mathcal O}_\varphi$
of $\varphi$ {\em i.\ e.} there exists a diffeomorphism $F$ of $M$ which sends
each leaf of ${\mathcal O}_\psi$ to a leaf of ${\mathcal O}_\varphi$.}
\end{Definition}

\begin{Definition} {\em An action $\varphi\in {\mathcal A}^r$ is said to be
{\em parameter rigid} if for any action $\psi\in {\mathcal A}^r$ such that
${\mathcal O}_\psi={\mathcal O}_\varphi$, there exist a diffeomorphism $F$ which
preserves leaves of these identical foliation and an automorphism
$A$ of $G$ such that (\ref{e1}) holds.}
\end{Definition}

The local version of parameter rigidity is not defined, simply because
we have no intermediate example. In Sect.\ 2, 3 and 10, we will
focus our attention
on the parameter rigidity for abelian or solvable Lie group actions.

\bigskip 

We will define the leafwise cohomology in Sect.\ 3 and discuss its close
relation with the parameter rigidity when the group $G$ is abelian.
After we introduced methods for the computation of the leafwise
cohomology in Sect.\ 4, Sect.\ 5, 6, 8 and 9 are devoted to
the computational results for various concrete foliations.

In Sect.\ 10, we raise  examples of parameter rigid solvable group actions.
The final Sect.\ 11 is devoted to the relation of the leafwise
cohomology to the problem of the existence of Riemannian metric
for which all the leaves are minimal surfaces.

\section{Regidity of flows} 

An $\R$-action $\varphi$ is what is usually called a flow, and is associated
with the vector field $X$ given by 
$$
X_x=\frac{d}{dt}\vert_{t=0}\varphi^t(x).
$$
An $\R$-action $\varphi$ is locally free if and only if $X$ is nonsingular.
As we have mentioned in Sect.\ 1, there is no $C^1$ locally rigid flow
unless the manifold is $S^1$. However there does exist a parameter rigid
flow, which we shall explain below.

Given a real number $\alpha\in\R$, define a Kroneker flow $\varphi_\alpha$
on the 2-torus $T^2$ by
$$
\varphi^t_\alpha(x,y)=(x+\alpha t, y+t).
$$
If the slope $\alpha$ is rational, then the flow is periodic {\em i.\ e.}
$\varphi^q$ is the identity for some $q>0$. All the orbits are closed and
their periods are the same.
The flow can never be parameter rigid, since one can change the
periods of closed orbits so that they are not identical.

If $\alpha$ is irrational, then
all the orbits are dense in $T^2$. It is also known that the flow is
uniquely ergodic {\em i.\ e.} the $\varphi_\alpha$-invariant probability is
unique. As for parameter rigidity, we have the following dichotomy.

\begin{Proposition} \label{p1}
The Kronecker flow $\varphi_\alpha$ is parameter rigid if and only
if the slope $\alpha$
is badly approximable.
\end{Proposition}

\begin{Definition} {\em A real number $\alpha$ is said to be {\em
badly approximable} if there exist $C>0$ and $\rho>0$ such that
$$
\Vert k\alpha\Vert_{S^1}>C\vert k\vert^{-\rho}, \ \ \ \forall k \in \Z\setminus
\{0\},
$$
where $\Vert \cdot \Vert_{S^1}$ denotes the distance to 0 of the projected 
image in $S^1=\R/\Z$.}
\end{Definition}

The proof of Proposition \ref{p1} in one direction is in order. Assume
$\alpha$ is badly approximable. Let $S^1$ be a circle in $T^2$ defined
by $y=0$. $S^1$ is a global cross section for $\varphi_\alpha$ and the first
return map is the rotation by $\alpha$, $R_\alpha$. Let $\psi$ be the flow
obtained by reparametrization of $\varphi_\alpha$, {\em i.\ e.} ${\mathcal O}
_\psi={\mathcal O}_{\varphi_\alpha}$. Let $f:S^1\rightarrow R$ be the first 
return time of $\psi$ for the cross section $S^1$; thus 
$\psi^{f(x)}(x)=R_\alpha(x)$.

\medskip
\noindent
\textsc{Claim:} {\em There exist $g\in C^\infty(S^1)$ and $c\in\R$ such that
$f=g\circ R_\alpha-g+c$.}

\medskip
Let us show first that Claim is sufficient to show the parameter rigidity
of $\varphi_\alpha$. For this, define a global cross section $C$ to the flow
$\psi$ by
$$
C=\{ \psi^{-g(x)}(x)\mid x\in S^1\}.
$$
Then we have
$$
\psi^c(\psi^{-g(x)}(x))=\psi^{-g(R_\alpha(x))}(\psi^{f(x)}(x))
=\psi^{-g(R_\alpha(x))}(R_\alpha(x)).
$$
This shows that the return time for the cross section $C$ is identically
equal to $c$. Now it is easy to construct a diffeomorphism $F$ conjugating
$\psi$ to $\varphi_\alpha$ up to time change $t\rightarrow ct$. ($F$ maps $C$ 
to $S^1$.)

\medskip
Let us turn to the proof of Claim. We shall show it for complex valued 
functions. A square integrable function $f\in L^2(S^1)$ can be expressed
as
$$
f=\sum_{k\in\Z}\hat f_k e^{2\pi i kt}, \ \ \ \hat f_k\in\C, \ \
\sum\vert \hat f_k\vert^2<\infty.
$$
Notice that $f\in C^\infty(S^1)$ if and only if for any $r>0$,
there exists $C_r>0$ such that $\vert \hat f_k\vert\vert k\vert^r<C_r$
for any $k\in Z$.
Now the equation $f=g\circ R_\alpha - g+c$ can be read off as
$$
\hat g_k=\hat f_k/(e^{2\pi i k \alpha}-1), \ \ \ c=\hat f_0,
$$
where $\hat g_k$ is the Fourier coefficients of $g$.
Since $\alpha$ is badly approximable and since
$$
\vert e^{2\pi ik\alpha}-1\vert\geq C_1\Vert k\alpha\Vert_{S^1}
$$
for some $C_1>0$, we have
$$
\vert \hat g_k\vert \vert k\vert^r\leq C_1^{-1}
 \frac{\vert \hat f_k\vert\vert k
\vert^r}{\Vert k\alpha\Vert_{S^1}}\leq C_2 \vert \hat f_k\vert \vert k\vert
^\rho\vert k \vert^r\leq C_3.
$$
This shows the smoothness of $g$, and the proof of Claim is complete.
We shall omit the proof of the other implication of Proposition \ref{p1}.

\bigskip

Badly approximability of the slope can be defined for Kronecker flows
on higher dimensional torus, and these flows are also shown to be
parameter rigid.
They constitute all the known examples. For related topics, see
\cite{AS,dS, LS}.

Here is a criterion of the parameter rigidity of a flow, which will be
useful for the study of dynamical properties of such flows.

\begin{Proposition} \label{p2}
The nonsingular flow $\varphi$ on $M$ 
given by a vector field $X$ is parameter rigid
if and only if for any $f\in C^\infty(M)$, there exist $g\in C^\infty(M)$
and $c\in\R$ such that $f=X(g)+c$.
\end{Proposition}

Let us show the only if part of the proposition. (The other part is just
by an analogous argument, using the integration instead
of the differentiation.) Choose $f\in C^\infty(M)$. We are going
to seek for $g$ and $c$ as in the proposition. It is no loss of generality
to assume that $f$ is positive. Let $\psi$ be the flow defined by
the vector field $(1/f)X$. Since $\psi$ and $\varphi$ has the same orbit
foliation, there is a function $\tau:\R\times M \rightarrow \R$ such that
$$
\varphi^t(x)=\psi^{\tau(t,x)}(x).
$$
Take $\frac{d}{dt}\vert_{t=0}$ of both sides, and we get
$$
\frac{d}{dt}\vert_{t=0}\tau(t,x)=f(x).
$$
Now by the parameter rigidity of $\varphi$, there is a diffeomorphism $F$ and
$c\in \R$ such that
$\psi^{ct}(F(x))=F(\varphi^t(x))$. The diffeomorphism $F$ preserves the orbits,
and can be written as $F(x)=\psi^{-g(x)}(x)$ for some $g\in C^\infty(M)$.
We have
\begin{gather*}
\psi^{ct-g(x)}(x)=\psi^{-g(\varphi^t(x))}(\varphi^t(x)), \ \ \rm{ and \ thus} \\
\psi^{g(\varphi^t(x))-g(x)+ct}(x)=\varphi^t(x).
\end{gather*}
That is,
$$
\tau(t,x)=g(\varphi^t(x))-g(x)+ct.
$$
Taking $\frac{d}{dt}\vert_{t=0}$, we obtain $f=X(g)+c$, as is desired.

\bigskip

An important implication of this observation is:

\medskip

\begin{Corollary} \label{c1}
A parameter rigid flow is uniquely ergodic, and it leaves smooth volume form
invariant.
\end{Corollary}

Notice that a flow $\varphi$ is uniquely ergodic  if and only if for any
$C^\infty$ map $f$, the orbit average
$$
\frac{1}{T}\int_0^Tf(\varphi^t(x)dt
$$
converges to a constant function, uniformly on $x\in M$.
Thus Propositon \ref{p2} immediately implies the unique ergodicity.

For the latter statement, let ${\rm vol}$ be a Riemannian
volume form on $M$. Take the Lie derivative:
$
L_X({\rm vol})=f\cdot{\rm vol}
$. By Proposition \ref{p2}, we have
$$
L_X({\rm vol})=(X(g)+c){\rm vol}.
$$
The integaration over $M$ shows that $c=0$. Thus we have
$$
X(e^{-g}\cdot{\rm vol})=0,
$$
{\em i.\ e.} $e^{-g}\cdot {\rm vol}$ is an invariant volume form.

\bigskip

However even among those flows which satisfy the necessary conditions of
the previous corollary, no examples of parameter rigid flows are found
except the linear flows on higher dimensional torus with badly approximable
slopes.
Here are two typical negative results in this direction.
Let $F_\lambda: T^2 \rightarrow T^2$ be a diffeomorphism defined by
$$
F_\lambda(x,y)=(x+y, y+\lambda),
$$
where $\lambda$ is a real number, and let $\varphi_\lambda$
 be the suspension flow.
For irrational $\lambda$, the flow $\varphi_\lambda$ satisfies the conditions of 
Corollary \ref{c1}.
However we have the following result found in \cite{K}.

\begin{Theorem} \label{t1}
For any real number $\lambda$, the flow $\varphi_\lambda$ is not parameter rigid.
\end{Theorem}

The other result is about the horocycle flows.
 Let us introduce them briefly.
Let $M$ be the quotient $\Gamma\setminus SL(2;\R)$ by a cocompact lattice
$\Gamma$. On $M$ the right action of the one parameter subgroup
$$
\{
\bmt 1 & t \\ 0 & 1\ema
\mid t\in\R
\}
$$
defines a flow, called {\em horocycle flow}.
Horocycle flows are known to satisfy
the conditions of Corollary \ref{c1} (\cite{F}).
The following theorem is due to L. Fluminio and G. Forni (\cite{FF}).

\begin{Theorem} \label{t2}
The horocycle
flow on a compact manifold $M$ is not parameter rigid.
\end{Theorem}

We will discuss these two theorems later in the next
section after the leafwise cohomology
of a foliation is introduced.

\smallskip

By the way, here is a very simple geometric proof of Theorem \ref{t2}
when the manifold $M$ is not a rational homology sphere. 

Let us take a basis of the Lie algebra $\mathfrak{sl}(2,\R)$ as follows.
$$
Y=\frac{1}{2}\bmt 1 &  \\ & -1 \ema, \ \ S=\bmt & 1 \\ & \ema, \ \
U=\bmt & \\ 1 & \ema
$$
These are left invariant vector fields on $SL(2;\R)$ and induces vector fields
on $M=\Gamma\setminus SL(2;\R)$. The horocycle flow $\varphi$ is the
one defined
by the vector field $S$.

Let us denote by $\eta$, $\sigma$ and $u$ 
the left invariant 1-forms on $SL(2;\R)$ which are
dual to $Y$, $S$ and $U$. They also induce 1-forms on $M$.

Since the first Betti number of $M$ is nonzero, there is a closed 1-form
$\omega$ such that the period map 
$$
[\omega]:\pi_1(M)\rightarrow \R
$$ 
is nontrivial and takes value in $\Z$. Let us write $\omega$ as
$$
\omega= f\sigma+*\eta+*u\ \ \ f,*\in C^\infty(M).
$$

Assume for contradiction that $f=S(g)+c$. Let $\hat M$ be the cyclic
covering of $M$ associated with the homomorphism $[\omega]$.
Then the lift $\hat \omega$ of $\omega$ is exact and the 
primitive $\hat h\in C^\infty(\hat M)$ is proper. Denote by $\hat \varphi$
the lift of the horocycle flow $\varphi$ to $\hat M$. For any $x\in M$,
any lift $\hat x$ of $x$, and any $T>0$, we have
\begin{gather}
\hat h(\hat\varphi^T(\hat x))-\hat h(\hat x)
=\int_{\hat\varphi^{[0,T]}(\hat x)}\hat \omega
=\int_{\varphi^{[0,T]}(x)}\omega
=\int_0^T f(\varphi^t(x))dt \notag \\
=\int_0^T (S(g)(\varphi^t(x))+c)dt
=g(\varphi^T(x))-g(x)+cT. \notag
\end{gather}

Now the lift of the horocycle flow $\hat \varphi$
has a dense orbit (Hedlund). Take $\hat x$
from a dense orbit. Assume $c\neq 0$.
Then there exists a sequence $T_i\rightarrow\pm\infty$ such that
$\hat \varphi^{T_i}(\hat x)\rightarrow \hat x$.  This is a contradiction
since the left hand side for $T_i$ tends to 0, while 
the right hand side tends to $\pm\infty$. Assume now $c=0$.
Then since the function $\hat h$ is proper, there exists $T'_j$ such 
that the left hand side for $T'_j$ tends to $\pm\infty$. Again a 
contradiction.

\smallskip

Here is a conjecture by A. Katok (\cite{K}).

\begin{Conjecture}
An arbitrary parameter rigid flow is smoothly conjugate to
a linear flow on the torus of badly approximable slope.
\end{Conjecture}

This conjecture is known to be true if the manifold has cup length
equal to the dimension, or if the manifold is  3-dimensional and
has nonvanishing Betti number.

\section{Leafwise cohomology and parameter rigidity}
In this section, we define the leafwise cohomology of foliations, and
discuss its  relationship with the
parameter rigidity  when the foliation is given
by a locally free $\R^p$-action.

Let $\mathcal F$ be a foliation on a manifold $N$, and denote by $T\mathcal F$
the tangent bundle of $\mathcal F$. A {\em leafwise $k$-form} $\omega$ 
is a $C^\infty$
cross section of the homomorphism bundle ${\rm Hom}(\bigwedge ^k
T{\mathcal F};\R)$. The space of leafwise $k$-forms is denoted by
$\Omega^k({\mathcal F})$.
 Thus given $k$ vector fields $X_1$, $X_2$, $\cdots$, $X_k$
tangent to $\mathcal F$ {\em  i.\ e.}  smooth cross sections of $T\mathcal F$,
and $\omega\in \Omega^k(\mathcal F)$, a smooth function $\omega(X_1, X_2,
\cdots. X_k)$ is defined.
The integrability condition for $T\mathcal F$ allows us to define
the {\em leafwise exterior derivative}
$$
d_\mathcal F: \Omega^k(\mathcal F)\rightarrow \Omega^{k+1}(\mathcal F)
$$
just as the usual exterior derivative. For example, for $\omega\in 
\Omega^1(\mathcal F)$ we define
$$
(d_\mathcal F \omega) (X_1,X_2)=X_1(\omega(X_2))-X_2(\omega(X_1))-\omega([X_1,X_2]).
$$

Then $(\Omega(\mathcal F),d_\mathcal F)$ constitutes a cochain complex,
whose cohomology is called the {\em leafwise cohomology} of $\mathcal F$,
denoted by $H^*(\mathcal F)$. This cochain complex is not elliptic, and
the leafwise cohomology can be infinite dimensional.

\begin{Example}
{\em The 0-dimensional leafwise cohomology $H^0(\mathcal F)$ coincides with
the vector space formed by {\em basic functions}
{\em i.\ e.} smooth functions on $M$, constant along
the leaves. Thus if $\mathcal F$ admits a dense leaf, then 
$H^0(\mathcal F)\cong\R$.}
\end{Example}

\begin{Example} \label{ex1}
{\em Let $N=L\times T$ and let $\mathcal F$ be the foliation on $N$ whose leaves
are $L\times\{ t\}$, $t\in T$. Then we have}
$$
H^k(\mathcal F)= H^k(L;C^\infty (T))=H^k(L;\R)\otimes C^\infty(T).
$$
\end{Example}

The cochain space $\Omega^k(\mathcal F)$ is equipped with the Whitney
$C^\infty$ topology, and the leafwise exterior derivative $d_\mathcal F$
is continuous. Thus the cocycle space ${\rm Ker}(d_\mathcal F)$ is
closed. But the coboundary space ${\rm Im}(d_\mathcal F)$ is not necessarily 
closed. The quotient of the cocycle space by the closure of the coboundary 
space is called the {\em reduced leafwise cohomology} and is denoted
by $\mathcal H^*(\mathcal F)$.

J. \'Alvarez L\'opez and G. Hector (\cite{AH})
 have given sufficient conditions for
the foliations to have infinite dimensional reduced cohomology,
and raised a lot of examples. Their examples are for the most part foliations
with dense leaves, definitely not like Example \ref{ex1}.
But in this talk, we are mainly interested in such a foliation for which
the leafwise cohomology is finite dimensional.

\bigskip

There are two ways, important for us,
to produce elements of the leafwise cohomology.
First notice that the restriction map
$$
r: {\rm Hom}(\bigwedge^*TM;\R)\rightarrow{\rm Hom}(\bigwedge^*T\mathcal F;\R)
$$
induces a cochain homomorphism (denoted by the same letter)
$$
 r:\Omega^*(M)\rightarrow \Omega^*(\mathcal F),
$$ 
where $\Omega^*(M)$ denotes
the de Rham complex of $M$. This induces a homomorphism of the 
cohomology groups
$$
r_*: H^*(M;\R) \rightarrow H^*(\mathcal F).
$$
The homomorphism $r_*$ is often nontrivial and yields elements
of $H^*(\mathcal F)$.

Secondly, suppose that the foliation $\mathcal F$ is given by
a locally free right action of a Lie group $G$. Then by the differentiation, we
get a Lie algebra homomorphism
$\iota: \mathfrak g \rightarrow \mathcal X^\infty(\mathcal F)$, where
$\mathfrak g=T_eG$ is identified with
the Lie algebra of the left invariant vector fields  on $G$,
and $\mathcal X^\infty(\mathcal F)$ the
Lie algebra of the vector fields of $M$ tangent to the foliation $\mathcal F$. 
Let $X_1,\cdots X_p$ be the basis of $\mathfrak g$. Then a leafwise
$k$-form $\omega$ is completely determined by $\omega(\iota X_{i_1},\cdots, 
\iota X_{i_k})$ for $1\leq i_1 < \cdots < i_k\leq p$, since for each point
$x\in M$ the tangent space $T_x(\mathcal F)$ is spanned by 
$(\iota X_1)_x, \cdots (\iota X_p)_x$. Thus given a left
invariant $k$-form 
$$\omega:\mathfrak g \times\cdots\times \mathfrak g \rightarrow \R,$$
a leafwise $k$-form $\iota \omega$ is defined by
$$
\iota\omega(\iota X_{i_1},\cdots,\iota X_{i_k})
=
\omega(X_{i_1},\cdots, X_{i_k}).
$$
This induces a homomorphism 
$$\iota_*: H^*(\mathfrak g)\rightarrow
H^*(\mathcal F).$$

\begin{Proposition}\label{p3}

(1) The homomorphism $\iota_*: H^1(\mathfrak g)\rightarrow
H^1(\mathcal F)$ is injective.

(2) If $G=\R^p$, then $\iota_*: H^i(\R^p)\rightarrow
H^i(\mathcal F)$  is injective for any $i\geq 0$.
\end{Proposition}

\begin{Remark} The cohomology of the abelian Lie algebra $\R^p$
is isomorphic to the exterior algebra of $\R^p$, and hence
to the cohomology of the $p$-torus:
$$
H^*(\R^p)\cong \bigwedge^*\R^p\cong H^*(T^p;\R).
$$
\end{Remark}

\noindent
Let us give a proof of Proposition \ref{p3}.
Let $\xi_1,\cdots, \xi_r$ be the closed left invariant 1-forms
whose classes form a basis of $H^1(\mathfrak g)$. They are of course
linearly independent in the dual space $\mathfrak g^*$ of $\mathfrak g$.
Therefore there exist elements $X_1,\cdots, X_r$ of $\mathfrak g$ such that
$\xi_i(X_j)=\delta_{ij}$.
Assume
$$
a_1\iota\xi_1+\cdots+a_r\iota\xi_r=d_\mathcal Ff
$$
for $a_i\in\R$ and $f\in C^\infty(M)$.
Let $\gamma_T:[0,T]\rightarrow M$ be an integral curve
of the vector field $\sum_ia_i\iota X_i$. Then we have
$$
f(\gamma_T(T))-f(\gamma_T(0))=\int_{\gamma_T}\sum_i a_i\iota\xi_i
=T\sum_i a_i^2.
$$
Since $T$ can be arbitrarily large and the left hand side is bounded,
all the coefficients $a_i$'s must vanish. This shows that the
classes of $\iota\xi_i$ are linearly independent in $H^1(\mathcal F)$.

The proof for the second part is basically similar and use the Stokes
theorem on an embedding of the rectangle $[0,T]\times\cdots\times [0,T]$. 

\bigskip

By Proposition \ref{p3}, if $\mathcal F$ is the orbit foliation of a locally 
free $\R^p$ action, then we have $\dim H^1(\mathcal F)\geq p$.
The following proposition relates the leafwise cohomology of the orbit
foliation to the parameter rigidity of the action.

\begin{Proposition} \label{p4}
A locally free effective $\R^p$-action is parameter rigid if and only if
its orbit foliation $\mathcal F$ satisfies $\dim H^1(\mathcal F)= p$.
\end{Proposition}

An action is said to be effective if the isotropy subgroup of some
point of the manifold is trivial.
The proof for $p=1$ is in order. Let $X$ be the vector field defining
a nonsingular flow $\varphi$. Define a leafwise 1-form $\omega$ by
$\omega(X)=1$. An arbitrary leafwise 1-form (which is always closed
by the dimension reason) is written as $f\omega$ for some $f\in C^\infty(M)$.
Given $g\in \Omega^0(\mathcal F)=C^\infty(M)$, notice that $(d_\mathcal F
g)(X)=X(g)$, and thus $d_\mathcal F g=X(g)\omega$. 
Therefore $\dim H^1(\mathcal F)=1$ if and only if for any
$f\omega$, there exist $g\in C^\infty(M)$ and $c\in\R$ such that
$$
f\omega=X(g)\omega +c\omega,
$$ 
which is equivalent to the parameter rigidity of the flow $\varphi$ by
Proposition \ref{p2}.
The proof for $p\geq 2$ is found in \cite{MM}.

\smallskip
Now let us return to the flows in Theorems \ref{t1} and \ref{t2}.
In fact what is proven respectively by A. Katok,
and L. Flaminio and G. Forni are much stronger than stated there.

\begin{Theorem}
Let $\mathcal F_\lambda$ be the orbit foliation of the flow
$\varphi_\lambda$ in Theorem \ref{t1}. Then the leafwise
cohomology $H^1(\mathcal F_\lambda)$ is infinite dimensional.
Moreover if $\lambda$ is badly approximable, then we have
${\mathcal H}^1(\mathcal F_\lambda)=H^1(\mathcal F_\lambda)$,
{\em i.\ e.} the space of the coboundaries is closed.
\end{Theorem}

\smallskip

\begin{Theorem}
Let $\mathcal F$ be the orbit foliation of the horocycle flow.
Then $H^1(\mathcal F)$ is infinite
dimensional and we have $H^1(\mathcal F)=\mathcal H^1(\mathcal F)$.
\end{Theorem}

\section{How to compute leafwise cohomology}

Here we will show some methods to compute the leafwise cohomology 
of a foliation. 

Let $\pi: M \rightarrow B$ be a {\em foliated bundle}
{\em i.\ e.} a fiber bundle equipped with a foliation  $\mathcal F$ 
on $M$ which is transverse to each fiber. Thus
the dimension of the leaves of $\mathcal F$ is greater than or
equal to the dimension of the base $B$. 

An important class of foliated bundes is obtainded by a construction called
suspension.
Let $(F,\mathcal G)$ be a foliated manifold.
Denote by ${\rm Diff}(F, \mathcal G)$ the group of diffeomorphisms
of $F$ which leaves the foliation $\mathcal G$ invariant.
Let $h:\Gamma\rightarrow {\rm Diff}(F, \mathcal G)$ be a homomorphism
from a group $\Gamma$. Let $B$ be a manifold such that $\pi_1(B)
\cong \Gamma$. Then a foliated manifold $(M,\mathcal F)$ called the
{\em suspension} of $h$ is defined as follows.
On the product $F\times\widetilde B$, where $\widetilde B$ stands
for the universal covering of $B$, there is defined a foliation
$\widetilde{\mathcal F}$ whose leaves are $L\times \widetilde B$, where
$L$ is a leaf of $\mathcal G$. Consider the diagonal action of 
$\Gamma$ on $F \times\widetilde B$, on the first factor through
$h$ and on the second by deck transformation.
This action clearly preserves the foliation $\widetilde{\mathcal F}$
and as the quotient we have a foliated manifold $(M, \mathcal F)$.

The following fundamental result is due to A. El Kacimi and A. Tihami
(\cite{ET}).

\begin{Theorem} \label{t3}
For the suspension foliation, there is a spectral sequence such that 
$$
E^{p,q}_2=H^p(B; H^q(\mathcal G))
$$
which converges to $H^{p+q}(\mathcal F)$, where $H^q(\mathcal G)$
is a system of local  coefficents on which $\Gamma=\pi_1(B)$ acts
through the homomorphism $h$.
\end{Theorem}

This shows for example that if 
$H^q(\mathcal G)$ is finite dimensional for any $0\leq q \leq n$,
then the leafwise cohomology $H^n(\mathcal F)$ is finite dimensional.

\smallskip
Let us consider a special case
where the foliation $\mathcal G$ on the fiber is a point foliation.
In this case we call $\mathcal F$ the {\em suspension of a point 
foliation} by a homomorphism
$h:\Gamma\rightarrow {\rm Diff}(F)$.
We have $H^k(\mathcal G)=0$ unless $k=0$. Thus the spectral sequence
of Theorem \ref{t3} collapses and we have:

\begin {Corollary} \label{c2}
If $\mathcal F$ is the suspension of a point foliation by a homomorphism
$h:\Gamma=\pi_1(B)\rightarrow {\rm Diff}(F)$
and if $\widetilde B$ is contractible, then we have:
$$
H^p(\mathcal F)=H^p(\Gamma; C^\infty(F)),
$$
for any $p\geq 0$, 
where $\Gamma$ acts on $C^\infty(F)$ through the homomorphism $h:
\Gamma\rightarrow {\rm Diff}(F)$.
\end{Corollary}

\bigskip

Here is another way for computing the leafwise cohomology, completely different
from the above, due to A. Haefliger (\cite{H}).

\begin{Theorem} \label{t5}
Let $\mathcal F$ be the  suspension of a point foliation by a homomorphism 
$h:\Gamma\rightarrow {\rm Diff}(F)$. Then there is an isomorphism
$$
\int_{\mathcal F}:H^{{\rm dim}\mathcal F}(\mathcal F) \rightarrow
H_0(\Gamma, C^\infty(F)).
$$
\end{Theorem}

The group $H_0(\Gamma, C^\infty(F))$ is by definition the quotient space
of $C^\infty(F)$ by the subspace spanned by the elements $f\circ h(\gamma)
-f$ for $f\in C^\infty(F)$ and $\gamma\in \Gamma$.

In sections below, we shall show  computational results for some
concrete foliations.

\section{Linear foliations}

Let us consider the suspension foliation of a point foliation,
where the group $\Gamma$ is $\Z^p$, the fiber $F$ is $T^q$,
and the homomorphism $h$  is given by
$$
h:\Z^p\rightarrow T^q \subset {\rm Diff}(T^q),
$$
 where $T^q$ is to be the group
of the translations of $T^q$.
Let $B$ be a real valued
$p\times q$ matrix. 

\begin{Definition}
{\em The matrix $B$ is called {\em badly approximable} if
there exist $C>0$ and $\rho>0$ such that
for any $k\in \Z^q\setminus\{0\}$, we have}
$$
\Vert Bk \Vert_{T^p}\geq C\vert k \vert^{-\rho}.
$$
\end{Definition}

Given a matrix $B$, a homomorphism 
$$h_B:\Z^p\rightarrow {\rm Diff}(T^q)$$
is defined by
$$
h_B(n)(x)=x+nB.
$$

As the suspension of $h_B$ over the base manifold $T^p$,
we get a linear foliation $\mathcal F_B$ on $T^{p+q}$.
The matrix $B$ is called the {\em slope matrix} of the foliation
$\mathcal F_B$.

\begin{Theorem}{\bf \cite{AS}} \label{t4}
If the slope matrix $B$ is badly approximable, then the homomorphism
$$
\iota_*: H^i(\R^p)\rightarrow H^i(\mathcal F_B)
$$
is an isomorphism for any $i\geq 0$.
\end{Theorem}

Recalling the definition of the homomorphism $\iota_*$, we get
the following corollary, which will be useful in Sect.\ 11.

\begin{Corollary} \label{c3}
Under the same assumption as above, the homomorphism
$$
r_*: H^i(T^{p+q})\rightarrow H^i(\mathcal F_B)
$$ is a surjection for any $i\geq 0$.
\end{Corollary}

As is easily shown the foliation $\mathcal F_B$ is the orbit
foliation of an $\R^p$ action. Since $H^1(\mathcal F_B)\cong \R^p$,
this action is parameter rigid if the action is effective.

In fact a bit wider class of $\R^p$-actions for which $\iota_*$
is an isomorphism are found in \cite{Lu}.

\section{Orbit foliation of transversely hyperbolic $\R^p$ actions}

Abundant examples of $C^1$ locally rigid (hence parameter rigid)
$\R^p$-actions are presented by A. Katok and others (\cite{K,
HK, KL,KS1,KS2}).
Here let us mention only one type among them.
Let us denote by $Aff_+(T^{p+1})$ the group of the orientation
preserving affine transformations of $T^{p+1}$, with respect to the
standard affine structure of $T^{p+1}$.
Let us consider a homomorphism 
$$
h: \Z^p\rightarrow Aff_+(T^{p+1})
$$
with the following property;

\smallskip

\noindent
(*) {\em The derivative $h_*:\Z^p\rightarrow SL(p+1;\Z)$
is injective and the image is generated by hyperbolic elements.}

\smallskip

In this case, the image $h_*(\Z^p)$ is shown to be 
simultaneously diagonalizable
in $SL(p+1;\R)$.

\begin{Theorem} {\bf \cite{KL}}
If $p\geq 2$, the 
cohomology $H^1(\Z^p; C^\infty(T^{p+1}))$ associated to $h$ above
is isomorphic to $\Z^p$.
\end{Theorem}

The suspension foliation $\mathcal F_h$
of $h$ over $T^p$ is the orbit foliation of an
$\R^p$ action $\varphi_h$. The above theorem implies that $H^1(\mathcal F_h)
\cong\R^p$, via Corollary \ref{c2}. That is, by Proposition
\ref{p4}, this action
is parameter rigid. Moreover it is shown to be $C^1$ locally rigid
(\cite{KL}).

\medskip

As for the higher cohomology group, however, this foliation does not
exhibit such rigid property as the foliation $\mathcal F_B$ does in Theorem
\ref{t4}.

\begin{Proposition} \label{p5}
For $p\geq 1$, the top dimensional cohomology $H^p(\mathcal F_h)$ 
of the foliation
$\mathcal F_h$ is  
infinite dimensional and we have $\mathcal H^p(\mathcal F_h)
=H^p(\mathcal F_h)$.
\end{Proposition}

The proof uses Theorem \ref{t5}.
For $p=1$ this result is classical.
For the related topics for Anosov flows, see \cite{L,GK,dLMM}.
\bigskip

Let us state a rigidity result about hyperbolic $\R^p$ actions,
which is global in nature.
An $\R^p$ action $\varphi$ on a closed $(2p+1)$-dimensional manifold $M$ is called
{\em split hyperbolic} if there is a continuous splitting
of the tangent bundle 
$$
TM=T\R^p\oplus E_1\oplus\cdots\oplus E_{p+1},
$$
where $T\R^p$ is the tangent bundle of the orbit foliation and each
$E_i$ is a 1-dimensional subbundle invariant under the differential
of the action. Furtheremore we assume that there exist $\xi_1,\cdots,
\xi_{p+1}\in\R^p$ such that the flow $\varphi(\cdot,t\xi_i)$ is expanding along 
$E_i$ and contracting along $E_j$ ($j\neq i$).

\begin{Theorem} {\bf \cite{M}}
A split hyperbolic $\R^p$ action is $C^\infty$ conjugate to the suspension
action $\varphi_h$ of some affine representation
$$
h:\Z^p\rightarrow Aff_+(T^{p+1})
$$ 
with the property (*), up to an automorphism of $\R^p$.
\end{Theorem}

\section{K\"unneth formula}

We shall interrupt the computation of the leafwise cohomology and
state a general theorem,
which has an interesting application for the parameter rigidity
of $\R^n$ actions.

A. El Kacimi and A. Tihami (\cite{ET}) followed the arguments of R. Bott and
L. W. Tu (\cite{BT}) in the framework of leafwise forms, to obtain {\em e.\ g.}
the Mayer-Vietoris theorem or the spectral sequence theorem
(Theorem \ref{t3}). One can further pursue this line to
obtain the following K\"unneth formula for the leafwise
cohomology.
Let $\mathcal F$ (resp.\ $\mathcal G$) be a foliation on
a manifold $M$ (resp.\ $N$). Then on the product manifold
$M\times N$, products of leaves of $\mathcal F$ and $\mathcal G$
form a foliation called the {\em product foliation}, denoted
by $\mathcal F\times \mathcal G$.

\begin{Theorem} \label{12}
Assume $\dim H^j(\mathcal G)<\infty$ for $0\leq j \leq k$. Then we
have an isomorphism:
$$
\sum_{i+j=k}H^i(\mathcal F)\otimes H^j(\mathcal G)\cong
H^k(\mathcal F\times \mathcal G).
$$
\end{Theorem}

Here is an immediate corollary concerning $\R^n$ actions..
Given two actions
$$
\varphi:M\times\R^p\rightarrow M \ \ {\rm and} \ \
\psi:N\times\R^q\rightarrow N,
$$
an $\R^{p+q}$ action $\varphi\times\psi$ on $M\times N$ called the
{\em product action} is defined by
$$
(\varphi\times\psi)((x,y),(s,t))=(\varphi(x,s),\psi(y,t)),
$$ 
where $x\in M$, $y\in N$, $s\in \R^p$ and $t\in\R^q$.

\begin{Theorem} Let $M$ and $N$ be closed manifolds.
Assume that $\varphi$ and $\psi$ are parameter rigid actions of
connected abelian Lie groups and that
the product action $\varphi\times\psi$ is effective.
Then $\varphi\times\psi$ is parameter rigid.
\end{Theorem}

For example the products of the linear actions of Sect.\ 5 and
the transversely hyperbolic  actions of Sect.\ 6 can be
parameter rigid.
It seems that this theorem is difficult to prove without resorting
to an algebraic topological argument.

Now in the next two sections, we will be back to the computation
of the leafwise cohomology.

\section{The weak stable foliation of the suspension Anosov flow}

Let $\mathcal W^s$ be the weak stable foliation  
of the suspension flow
of a hyperbolic toral automorphism $A: T^n\rightarrow T^n$.
Let $\mathcal V^s$ be the linear foliation on $T^n$ whose leaves
are parallel translations of the stable eigenspace $E^s$.
The automorphism $A$ leaves $\mathcal V^s$ invariant,
and as is well known the slope matrix of
the linear foliation $\mathcal V^s$ is badly approximable.
By Theorem \ref{t4}
we have 
$$
H^*(\mathcal V^s)\cong H^*(\R^p)\cong\bigwedge^*\R^p,
$$
 where $p$ denotes the dimension
of the foliation $\mathcal V^s$. The foliation $\mathcal W^s$ in
question is the suspension of $(T^n, \mathcal V^s)$ by a homomorphism
$h:\Z\rightarrow {\rm Diff}(T^n,\mathcal V^s)$ given by $h(1)=A$.
In the case of foliated bundles over $S^1$,
Theorem \ref{t3} gives rise to the Wang exact sequence;
$$
 H^{k-1}(\mathcal V^s) \rightarrow H^{k-1}(\mathcal V^s)
 \rightarrow H^k(\mathcal W^s)
\rightarrow H^k(\mathcal V^s) \rightarrow H^k(\mathcal V^s).
$$
The first and the last arrows are $A^*-{\rm Id}$.
Notice that $A^*:H^i(\mathcal V^s)\rightarrow H^i(\mathcal V^s)$ 
is the $i$-th exterior product of $A\vert_{E^s}$, under the identification
given by $\iota_*$: $\bigwedge^i E^s\cong H^i(\R^p)
\cong H^i(\mathcal V^s)$. Thus we get:

\begin{Proposition} \label{pa1}
$$
H^i(\mathcal W^s)=\R \ {\rm for} \  i=0,1, \ \  {\rm and}  \ \
H^i(\mathcal W^s)=0 \ {\rm for} \ i\geq 2.
$$
\end{Proposition}
We shall return to this result in Sect.\ 11.

\section{Stable foliations of geodesic flows}

There are foliations for which neither Corollary \ref{c2}
nor Theorem \ref{t5} apply but still we can compute the leafwise
cohomology by a geometric method.

Let $M=\Gamma\setminus {\rm SL}(2;\R)$ where $\Gamma$ is a cocompact
lattice and let $\mathcal F^s$ be the foliation defined by
the action of the subgroup
$$
\{ \bmt e^{t/2} & x \\ 0 & e^{-t/2}\ema \mid t,x\in\R \}.
$$

This subgroup is isormorphic to the Lie group $GA$ of the orientation
preserving affine action of the real line and we have
$H^1(\mathfrak g\mathfrak a)\cong\R$.
The foliation $\mathcal F^s$ is the weak stable foliation of
an Anosov flow.

\begin{Theorem} \label{t6}
The homomorphism
$$
(r_*,\iota_*): H^1(M;\R)\oplus\R\rightarrow H^1(\mathcal F^s)
$$
is an isomorphism.
\end{Theorem} 

Since there are many \lq\lq leafwise 1-currents"
given by closed orbits of the Anosov flow,
it is not difficult to show that $(r_*,\iota_*)$ is injective.
To show the surjectivity, we use the following lemma.

\begin{Lemma} 
Given a $C^\infty$ function $k$, the equation
\begin{equation} \label{e2}
k=Y(g)+g
\end{equation}
has a $C^\infty$ solution $g$ if and only if there exists a $C^\infty$
function $l$ fulfilling
$$
S(k)=Y(l).
$$
\end{Lemma}

Here $Y$, $S$ and $U$ are the vector fields on $M$ induced by
the action of the elements of the Lie algebra ${\mathfrak sl}(2;\R)$;
$$
Y=\frac{1}{2}\bmt 1 & 0 \\ 0 & 1\ema, \ \ 
S=\bmt 0 & 1 \\ 0 & 0\ema, \ \
U=\bmt 0 & 0 \\ 1 & 0 \ema.
$$

The integration along the orbit of the $Y$-flow from $t=-\infty$
gives a continuous solution $g$ of (\ref{e2}), which is $C^\infty$
 along $Y$ and
$U$ directions. The condition $S(k)=Y(l)$ transforms
the difference of $g$ between
nearby two points on an $S$-leaf 
to the sum of the difference of the images of these points by time $-T$ of
the $Y$-flow and the integration of $k$ along this long $S$-orbit.
Then the unique ergodicity of $S$-flow (\cite{F}) ensures that $g$ is
$C^\infty$ along the $S$-direction.

\medskip

Here is also a result about 2-dimensional cohomology, which is a
translation  of a result of A. Haefliger found in \cite{H}. 
An application of Theorem \ref{t5}.

\begin{Theorem} \label{t7}
The homomorphism
$r_*: H^2(M;\R)\rightarrow \mathcal H^2(\mathcal F^s)$
is an isomorphism.
\end{Theorem}

Unfortunately the above theorem is only about the reduced
leafwise cohomology $\mathcal H^2(\mathcal F^s)$. The honest
leafwise cohomology $H^2(\mathcal F^s)$ is not known.
We will be back to this problem in Sect.\ 11.

\section{Parameter rigidity of certain solvable Lie group actions}

First of all let us state a criterion for the parameter rigidity of 
the action of a general Lie group, which 
takes form of a nonlinear equation if the Lie group is nonabelian.

Let $\varphi:M\times G\rightarrow M$ be a locally free action of a 
connected and simply connected Lie group $G$ with Lie algebra $\mathfrak
g$. Denote by $\mathcal F$ the orbit foliation.
As is explained in Sect.\ 3, the vector field $\iota(X) \in T\mathcal F$
is defined for each $X\in\mathfrak g$. 
A $\mathfrak g$-valued leafwise 1-form $\omega_0$
associted with the action $\varphi$
is defined by
$$
\omega_0(\iota(X))=X.
$$
It satisfies
$$
d_{\mathcal F}\omega_0+[\omega_0\wedge\omega_0]=0,
$$
where $[\omega_0\wedge\omega_0]$ is a leafwise 2-form defined by
$$
[\omega_0\wedge\omega_0](X,Y)=[\omega_0(X),\omega_0(Y)].
$$

\begin{Proposition} \label{p6} Assume the action $\varphi$ is effective.
The action $\varphi$ is parameter rigid if for any $\mathfrak g$-valued
leafwise 1-form $\omega$ such that
$$
d_{\mathcal F}\omega+[\omega\wedge\omega]=0,
$$
there exist an endomorphism $\Phi$ of $G$ and a smooth map $b:
M\rightarrow G$ such that
$$
\omega=b^{-1}\Phi_*(\omega_0)b+b^{-1}db.
$$
\end{Proposition}

As an application let us consider the foliation of Sect.\ 8.
Computation based on the result of Theorem \ref{t6} leads to
an alternative proof of the following theorem by E. Ghys (\cite{G}).

\begin{Theorem} \label{11}
The action of the Lie group $GA$ in Sect.\ 8 is parameter rigid,
if the manifold $M$ is a rational homology sphere.
\end{Theorem}

E. Ghys also showed that if the manifold $M$ has nonvanishing
first Betti number, then the acion of $GA$ is parameter rigid
among those actions which preserve the volume form.
Our method using the leafwise cohomology just yields the same
conclusion. Nothing more, except we can show that the sufficient
condition of Proposition \ref{p6} are not
satisfied in this case. So the question of the parameter rigidity
of the $GA$-actions for 3-manifolds with nonvanishing first Betti
number is still open.

\medskip
Let us raise other examples. Let $A$ be a hyperbolic automorphism of the
torus $T^n$ ($n\geq 2$) and let $\mathcal W^s$ be the weak stable foliation
of the suspension flow. By Proposition \ref{pa1}, we have
$H^1(\mathcal W^s)\cong\R$. If the matrix $A$ has no negative
eigenvalues of absolute value smaller than 1, then the foliation $\mathcal W^s$
is the orbit foliation of a locally free action of a two step solvable group
$G$.

\begin{Theorem} {\bf \cite{MM}}
If the characteristic polynominal of $A$ is irreducible over ${\mathbb Q}$,
then the $G$ action is parameter rigid.
\end{Theorem}

Together with the rigidity result of the orbit foliations established
by A. El Kacimi and M. Nicolau (\cite{EN}), we obtain the following.

\begin{Theorem}
If furthermore all the eigenvalues of $A$ are positive, then the $G$
action is $C^\infty$ rigid.
\end{Theorem}

\section{Minimizable foliation}

A foliation $\mathcal F$ on a manifold $M$ is called {\em
minimizable} if there is a Riemannian metric on $M$ for which
all the leaves of $\mathcal F$ are minimal surfaces.

\begin{Proposition}
A foliation $\mathcal F$ on the total space $E$ of a fiber bundle
$\pi:E\rightarrow B$ transverse to each fiber and of complementary
dimension is minimizable.
\end{Proposition}

In fact a Riemannian metric on $B$ can be lifted to a leafwise
Riemannian metric of $\mathcal F$ and extended to a Riemannian metric
of $E$ in such a way that each fiber is orthogonal to each leaf.
Clearly any leaf is totally geodesic for this metric.

Thus most of the examples of foliations so far listed are minimizable.
However there is an exception.

\begin{Proposition}
Horocycle flows on  closed 3-manifolds are not minimizable.
\end{Proposition}

In fact a foliation having a positive holonomy invariant measure which
is the boundary of an invariant 1-current is not minimizable (\cite{S}).

Here is a criterion for the minimizability due to
H. Rummler and D. Sullivan (\cite{R,S}). Assume all the foliations
in this section are oriented.

\begin{Theorem} \label{13}
Let $g_0$ be a leafwise Riemannian metric of a $p$-dimensional
 foliation $\mathcal F$.
It extends to a Riemannian metric of the whole manifold $M$
for which all the leaves are minimal surfaces if and only if
the leafwise volume form $\omega_0$ induced by $g_0$
extends to a  form $\omega$ of $M$ such that
$$
d\omega(\xi_1,\cdots,\xi_{p+1})=0
$$
whenever the first $p$ of the vector fields $\xi_i$'s are tangent to
the leaves.
\end{Theorem}

Now we call a foliation $\mathcal F$ {\em totally minimizable}
if any leafwise Riemannian metric extends to a Riemannian
metric for which all the leaves are
minimal. 

\begin{Corollary} \label{14}
A $p$-dimensional foliation $\mathcal F$ on $M$ is totally minimizable
if the homomorphism
$$
r_*:H^p(M;\R)\rightarrow H^p(\mathcal F)
$$
is surjective. If ${\rm codim}\mathcal F=1$, then the converse is also true.
\end{Corollary}
 
In fact if $r_*$ is surjective, then the leafwise volume form
$\omega_0$ of any leafwise Riemannian metric can be expressed as
$$
\omega_0=r(\omega)+d_{\mathcal F}(\eta_0),
$$
for a closed $p$-form $\omega$ and a leafwise $(p-1)$-form $\eta_0$.
But since $r:\Omega^{p-1}(M)\rightarrow\Omega^{p-1}(\mathcal F)$
is surjective, we have $\eta_0=r(\eta)$ for some $(p-1)$-form $\eta$.
It follows then that
$$
\omega_0=r(\omega+d\eta),
$$
for a closed $p$-form $\omega+d\eta$, whence the criterion of Theorem \ref{13}
is satisfied.

To show the converse statement for codimension one foliation, notice
first of all that since $\dim(M)=p+1$, the condition for
the $(p+1)$-form $d\omega$ in Theorem \ref{13} is 
equivalent to saying that $d\omega=0$.
Next notice that any leafwise $p$-form is the difference of
the leafwise volume forms of two leafwise Riemannian metrics.

\medskip

As applications, 
 Corollary \ref{c3} implies that the linear foliations on the torus 
with badly approximable slope
matrices are totally minimizable.
Also the foliation 
$\mathcal W^s$ of Sect.\ 8. is totally minimizable by Proposition \ref{pa1}.

\medskip
A foliation $\mathcal F$ is said to be {\em almost totally minimizable}
if any leafwise Riemannian metric is approximated arbitrarily
close in the $C^\infty$ topology by the restriction of  Riemannian
metrics for which all the leaves are
minimal. 

One example of such foliations are $p$-dimensional 
linear foliations $\mathcal F$
on the torus $T^n$
whose leaves are all dense. This is a consequence of the
surgectivity of the map $r_*:H^p(T^n)\rightarrow \mathcal H^p(\mathcal F)$.
This result is best possible if the slope matrix of the foliation
$\mathcal F$ is not badly approximable. (Otherwise $\mathcal F$ 
is totally minimizable, as we mentioned above.)

Likewise for the foliation $\mathcal F^s$ of Sect.\ 9, 
the homomorphism
$r_*: H^2(M;\R)\rightarrow \mathcal H^2(\mathcal F^s)$
is an isomorphism (Theorem \ref{t7}). This implies that the foliation
$\mathcal F^s$ is almost totally minimizable (\cite{H,HL}).
Unlike the previous example of linear flows, we do not know
if this is best possible or not.
It is known that the subspaces $\{Y(f)\vert f\in C^\infty(M)\}$ and
$\{S(f)\vert f\in C^\infty(M)\}$ are closed in $C^\infty(M)$
(\cite{GK} or Proposition \ref{p5}, and \cite{FF} or Theorem \ref{t2}).
The problem is equivalent to deciding if the sum of these
two subspaces is closed.

\end{document}